\documentclass[10pt,a4paper]{article}
\usepackage{amsmath}
\usepackage{latexsym}
\usepackage{theorem}
\usepackage[all]{xy}
\newtheorem{teorema}{Theorem}[section]
\newtheorem{definicion}[teorema]{Definition}
\include{amslatex}
\newtheorem{proposicion}[teorema]{Proposition}
\newtheorem{lema}[teorema]{Lemma}
\newtheorem{corolario}[teorema]{Corollary}
\textheight
24 cm \textwidth 14 cm \oddsidemargin 10 mm \topmargin -15 mm
{\theorembodyfont{\rmfamily}
\newtheorem{comentario}[teorema]{Remark}}
{\theorembodyfont{\rmfamily}
\newtheorem{ejemplo}[teorema]{Example}}

\numberwithin{equation}{section}

\include{amsmath}
\begin{document}
\begin{title}
{\large {\bf On Fermat's principle for causal curves in time oriented Finsler spacetimes}}
\end{title}
\maketitle

\begin{center}
Ricardo Gallego Torrom\'e\footnote{{Instituto de Matem\'atica e Estat\'istica - USP, S\~ao Paulo, Brazil.} email:
rgallegot@gmx.de. Financially supported by FAPESP, process 2010/11934-6},\,
Paolo Piccione\footnote{{Instituto de Matem\'atica e Estat\'istica - USP, S\~ao Paulo, Brazil}. Partially Sponsored by CNPq and Fapesp, Brazil},
\,Henrique Vit\'orio\footnote{{Instituto de Matem\'atica e Estat\'istica - USP, S\~ao Paulo, Brazil. Financially supported by CNPq, process 150124/2011-2}}
\end{center}

\abstract{\small{In this work, a version of {\it Fermat's principle} for causal curves with the same energy in time orientable Finsler spacetimes
is proved. We calculate the second
variation of the {\it time arrival functional} along a geodesic in terms of the index form associated with the Finsler spacetime Lagrangian.
 Then the character of the critical points of the time arrival functional is investigated and a
 Morse index theorem in the context of Finsler spacetime is presented.}}

\section{\small{Introduction}}
Finsler spacetimes appear in a natural way when modeling several physical phenomena. For instance, they are useful for the investigation of the propagation of light in locally anisotropic media (see for instance \cite{Perlick2000, Perlick06}), the investigation of multi-refringence \cite{SkakalaVisser, SkakalaVisser2} or as geometric models for classical point electrodynamics Randers spaces \cite{GallegoTorrome, Randers}. Recently, several results of {\it phenomenology of quantum gravity} have been related with Finsler spacetimes. It turns out that quantum gravity models generally predict modified dispersion relations at low energy.
 The investigation of such dispersion relations has shown that they can be associated
 with (in general not regular) {\it Finsler geometries} \cite{GirelliLiberatiSindoni}.
 Another example where Finsler spacetimes appear is in relation with
 the theory of {\it very special relativity} of Cohen and Glashow \cite{CohenGlashow}. It was shown by Gibbons et al. \cite{GibbonsGomisPope}
 that such spacetime geometries correspond to  Finsler spacetimes of Bogoslosvky's type \cite{Bogoslosvky}. Applications to cosmology has been found in \cite{KouretsisStathakopoulosStraviros}. It is also of relevance the relation with Lorentz violations models (see for instance \cite{Kosteleky}).

The theory of Finsler geometry of positive definite metrics admits a complete set of tools to be investigated
 in a similar way as for Riemannian geometry \cite{BCS}. In a less developed stage
 is the geometry  of {\it Finsler spacetimes}. One of the
  objectives of the present work is to show how some relevant methods can be
 {\it transported} from Lorentzian geometry to Finsler spacetime geometry in a natural way. In particular, we will concentrate on the extension of
 Fermat's principle to Finsler spacetimes and the associated variational theory.

In the framework of static solutions of the Einstein equations,
Fermat's principle appears in
the work of H. Weyl \cite{Weyl}. For light-like geodesics in a general spacetime, Fermat's principle was formulated mathematically by V. Perlick \cite{Perlick},
 as an attempt to formalize the theory of I. Kovner \cite{Kovner}. A complementary version of Fermat's principle for timelike curves in a time-oriented spacetime was investigated by \cite{GiannoniMasielloPiccione},
among others generalizations appearing in the literature.  For lightlike curves, there is already a version of Fermat's principle for Finsler spacetimes in the sense of J. Beem's \cite{Perlick06}. In this paper we present a  Fermat's principle for both lightlike and timelike Finsler in time oriented Finsler spacetimes. We also discuss
  the character of the critical points of the time arrival functional and we describe how to obtain by using standard methods borrow from Lorentzian geometry a Morse index theorem for the time arrival functional for timelike curves.

 Currently, there are several frameworks for Finsler spacetimes in the literature where notions
   of lightlike, timelike and spacelike curves and causal structure are available.
   The definition of Finsler spacetime that we took was introduced by
    J. Beem \cite{Beem1}. With such definition one can extend the results of Perlick \cite{Perlick} for causal curves (both timelike and lightlike)
   and obtain a formula for the second variation of the {\it time arrival functional}.

   Beem's framework does not contemplate all the Finsler spacetime models appearing in physical applications. For instance, a convenient way to describe the motion of a point charged particle in a external field is by a Randers metric \cite{Randers}.
However, it is still an open problem if Randers spaces of Lorentzian signature are naturally described in Beem's formalism. Also, some of the Finsler spacetime geometries associated with phenomenology of quantum gravity
  contain {\it singular sectors} in the tangent bundle of the base manifold, where the fundamental tensor is not enough regular. In order to deal with such models,
  one needs to consider weaker hypothesis in  Beem's definition, restricting the domain of definition of the Lagrangian $L$ to
  regions where the geometric objects (metrics, connections, etc) have enough regularity. We demonstrate that this natural modification of Beem's theory constitutes a convenient framework to investigate mathematical properties of
    Finsler spacetimes and in particular, to investigate Fermat's principle and related results.
\section{\small{Geometric framework}}
\subsection{\small{Finsler spacetimes}}
Following  J. Beem \cite{Beem1},
we introduce the basic notation and fundamental notions of Finsler spacetimes.
Let $M$ be a differentiable manifold and $TM$ the tangent bundle of $M$. Local
coordinates $(U,x)$ on $M$ induce local natural coordinates $(TU, x, y)$ on $TM$. The slit tangent bundle is $N=TM\setminus\{0\}$,
where $0$ is the zero section of $TM$.
\begin{definicion}
A Finsler spacetime is a pair $(M,L)$ where
\begin{enumerate}
\item $M$ is an $n$-dimensional real, second countable, Hausdorff $C^{\infty}$-manifold.
\item $L:N\longrightarrow R$ is a real smooth function such that
\begin{enumerate}
\item $L(x,\cdot)$ is positive homogeneous of degree two in the variable $y$,
\begin{align}
L(x,ky)=\,k^2\,L(x,y),\quad \forall\, k\in ]0,\infty[,
\label{homogeneouscondition}
\end{align}
\item The {\it vertical Hessian}
\begin{align}
g_{ij}(x,y)=\,\frac{\partial^2\,L(x,y)}{\partial y^i\,\partial y^j}
\label{nondegeracy-signature}
\end{align}
is non-degenerate and with signature $(-,+,...,+)$ for all $(x,y)\in\, N$.
\end{enumerate}
\end{enumerate}
\label{Finslerspacetime}
\end{definicion}
Direct consequences of this definition and Euler's theorem for positive homogeneous functions are the following relations,
\begin{align}
\frac{\partial L(x,y)}{\partial y^k}\,y^k=\,2\,L(x,y),\quad \frac{\partial L(x,y)}{\partial y^i}=\,g_{ij}(x,y)y^j,
\quad L(x,y)=\frac{1}{2}\,g_{ij}(x,y)y^iy^j.
\end{align}
Note that because homogeneity property on $y$, it is equivalent to provide the Lagrangian $L(x,y)$ or the corresponding fundamental tensor $g_{ij}(x,y)$.

There are other definitions of Finsler spacetime in the literature with its own merits.
 One of them goes back to Asanov \cite{Asa}.
His notion of Finsler spacetime is useful when dealing with timelike trajectories, and it can
be applied to investigate for instance timelike curves in {\it Randers type} Finsler spacetimes \cite{Randers}. However, it does not allow a covariant
notion of lightlike vectors and curves. A related theory of Finsler spacetimes is the notion {\it conic Finsler spaces}, developed by Javaloyes and Sanchez \cite{JavaloyesSanchez}. In such formulation, the timelike vectors at each point of the spacetime are defined on an open cone on each tangent space $T_xM$. It does not consider lightlike vectors and curves.
Other useful theory of Finsler spacetimes was introduced  by
Pfeifer and Wohlfarth \cite{PfeiferWohlfarth}. Their theory allows to consider some relevant Finsler spacetimes that are outside Beem's framework (for instance, bi-metric spaces). Still, Pfeifer-Wohlfarth theory does not contemplate lightlike curves in  Randers types of Finsler spacetimes.

The following is the notion of reversibility that we will consider,
\begin{definicion}
A Finsler spacetime $(M,L)$ is {\it reversible} iff
\begin{align}
L(x,y)=\,L(x,-y)
\label{reversibilityofL}
\end{align}
 for any $(x,y)\in \,N$. Otherwise, $(M,L)$ is a {\it non-reversible} Finsler spacetime.
 \label{definitionreversivilityofL}
\end{definicion}
This notion of reversibility is different from the one considered by Beem \cite{Beem1} and also different than the one considered in \cite{PfeiferWohlfarth}. Our definition of reversible metric is stronger than the corresponding notions of Beem and that the one considered in the theory of Pfeifer-Wohlfart.
 \subsection{\small{Elementary causality notions for Finsler spacetimes}}
The fundamental causal notions of a Finsler spacetime
$(M,L)$  is a natural generalization of the Lorentzian causal framework \cite{BEE}.
A vector field $X\in\,\Gamma TM$ is said to be timelike if $L(x, X(x))<0$ at all points $x\in\,M$ and a curve $\lambda:I\longrightarrow M$ is timelike
if the tangent vector field is timelike $L(\lambda(s),\dot{\lambda}(s))<0$. A vector field $X\in \,
 \Gamma$ is lightlike if $L(x,X(x))=0,\,\forall x\in\,M$; a curve is lightlike if its tangent
  vector field is lightlike. Similar notions hold for spacelike vector and curves.
A curve is causal if it is either timelike and has constant
  speed $g_{\dot{\lambda}}(\dot{\lambda},\dot{\lambda}):=L(\lambda,\dot{\lambda})=g_{ij}(\dot{\lambda},T)\dot{\lambda}^i\dot{\lambda}^j$ or if it is lightlike.

The following facts can be proved from the definition of Finsler spacetime,
\begin{enumerate}
\item The function $L(x,y)$ defines a positive definite, homogeneous function $F(x,y)=\,\sqrt{-L(x,y)}$ of degree one on the sub-bundle of timelike vectors
    \begin{align*}
    T^+M:=\,\{(x,y)\,\in TM,y\,\,\in T_x M\,s.t.\,\,L(x,y)<0\}.
\end{align*}
\item Each connected component of $T^+_xM$ is an open convex cone \cite{Perlick06}.
\end{enumerate}
The {\it Finsler function} $F(x,y)=\,\sqrt{-L(x,y)}$  defined on $T^+_xM$
 determines a Finsler spacetime in the sense of Asanov. The function $F$ can be extended in a non-smooth way to the whole bundle $TM$;
for a spacelike vector it is defined by $F(x,y)=\,\sqrt{L(x,y)}$ for spacelike vectors $y$ such that $L(x,y)$ and $zero$ for vectors $y$ in the null
cone,
\begin{align}
\mathcal{C}:=\,\{(x,y)\,\in TM,\,y\,\in T_x M\,s.t.\,\,L(x,y)=0\}.
\end{align}

A time  orientation is a smooth timelike
vector field $ T\in\Gamma TM$,
\begin{align}
g_T(T,T)<0.
\end{align}
\begin{definicion}
A timelike tangent vector $v\in T_xM$ is future pointed respect to $T$ if $g_v(v,T)<0$. It is easy to see that the set of future pointed vectors respect to $T$ form an open sub-set of the space of timelike vectors.
A timelike curve $\lambda:I\longrightarrow M$ is {\it future pointed} respect to $T$ if its tangent vector field is future pointed,
\begin{align*}
g_{\dot{\lambda}}(\dot{\lambda},T)=\,g_{ij}(\lambda(t),\dot{\lambda}(t))\,T^i\,\dot{\lambda}^j<\,0.
\end{align*}
 \label{definicionoftimelikeorientatio}
  \end{definicion}
  Similar notions hold for lightlike vectors and curves and for past pointed vectors.
  Note that for general Finsler spacetime $(M,L)$, given a tangent vector $v$ in the connected cone component where $T$ is contained, one does not expect that is future pointed respect to $T$. However, for close enough curves $\tilde{\gamma}:I\to M$ to $\gamma:I\to M$, tangent vectors $\tilde{v}$ to $v$, it holds by continuity of the function
  \begin{align*}
  g_{\cdot}(T,\cdot):\,T^+_xM \to R,\quad \tilde{v}\mapsto g_{\tilde{v}}(\tilde{v}, T)
  \end{align*}
  that if $g_v(v,T)<0$, then $g_{\tilde{v}}(\tilde{v}, T)<0$. This fact is needed when doing variation calculations.

  An {\it observer} is described by a future pointed timelike curve $\gamma:[a,b]\to M$;if there is a time orientation $T$, an observer is future pointed iff $g_{\gamma'}(\gamma',T)<0$.
  where $\gamma'$ is the tangent vector associated with the observer $\gamma$.

  In general, an observer is not necessarily described by integral curves of the time orientation $T$. Also, given a future pointed vector $w\in \, T^+_x M$ it can happen that for two observers $\gamma$ and $\tilde{\gamma}$, one has that $g_{w}(w,\gamma')<0$ but for the second observer one can have that $g_{w}(w,\tilde{\gamma})> 0$. This fact implies that one cannot decide the future pointed character of the vector $w$ from observations (that is, from measurements made by a given observer). Since the character of the sign of $g_w(w,T)$ is not controlled by the sign of the observable $g_w(w,\gamma')$, one needs to assume the sign of $g_w(w,T)$. Also, note the following fact: in general, one has that
  \begin{align*}
  g_w(w,-Z)=\,-g_w(w,Z).
  \end{align*}
  In particular, given a time orientation $T$, one has that $g_w(w,-T)=\,-g_w(w,T)$ and $g_{\gamma'}(\gamma',T)=-\,g_{\gamma'}(\gamma',T)$.

  In contrast, one has in general that
  \begin{align}
  g_{-w}(-w,Z)\neq \,- g_{w}(w,Z)
  \label{generalparityasymmetry}
  \end{align}
  for any $Z\in \Gamma N$. In particular, one has the property that for a time orientation $T$ and an observer $\gamma$, one has in general that
  \begin{align}
  g_{-\gamma'}(-\gamma',T)\neq \,- g_{\gamma'}(\gamma',T).
  \end{align}
 Also surprising, when $L$ is not reversible, if $v$ is a causal vector, then $-v$ is not necessarily causal. An example of non-reversible Finsler spacetime is provided by Ruth solution of a Finsler generalization of Einstein equations discussed bellow (see the space defined by equation \eqref{rutzexample}).

One way to avoid these puzzling consequences is to consider reversible metrics:
  \begin{proposicion}
  If the Finsler spacetime $(M,L)$ is reversible, then
  \begin{align}
  g_{-w}(-w,Z)=\,- g_{w}(w,Z).
  \label{symmetrycondition}
   \end{align}
   In particular, for any time orientation $T$ and observer $\gamma$, one has that
   \begin{align*}
   g_{-\gamma'}(-\gamma',T)= \,- g_{\gamma'}(\gamma',T).
   \end{align*}
  \end{proposicion}
  {\bf Proof}. That $(M,L)$ is reversible  means that $L(x,y)=\,L(x,-y)$. Therefore, for any vector fields $w$ and $Z$, one has that
  \begin{align*}
  g_{-w}(-w,Z)=\,g_{w}(-w,Z)=\,-g_w(w,Z).
  \end{align*}
  The second equality follows directly from \eqref{symmetrycondition}.\hfill$\Box$

  Note that the symmetry condition \eqref{symmetrycondition} hold for Lorentzian spacetimes, which makes
  the condition for $(M,L)$ being reversible a natural condition for some physical models. However, we observe that reversibility of $L$ is indeed not required for the proofs of the main results of this paper. Therefore, we will consider Finsler spacetimes that could also be non-reversible. Let us note that Randers spacetimes were introduced as a model where irreversibility in evolution was contained in the geometry of the spacetime \cite{Randers}. Although this is not necessarily our position, we should be alert that any un-necessary restriction could private us of a framework to explore general physical models.

  Despite these subtleties that surround  the notion of time orientation and non-reversibility in Finsler spacetimes, the notion of future pointed vector is a geometric notion in the sense that does not depend on the observer $\gamma$. It depends on the vector field $T$ and the vector $w$ only. In order to make sense of this notion, the vector field $T$ needs to be fixed and should be measurable for any physical observer $\gamma$. In particular, any future pointed observer should agree on the criteria $g_{\gamma'}(\gamma',T)<0$. By the discussion above, for each $x\in T_xM$, for each $x\in M$, this select the curves from the open convex component of the timelike vectors at $x$ containing $T(x)$.

  Once $T$ is fixed, the future pointed observers are well defined, as well as the past pointed observers. However, for a generic Finsler spacetime $(M,L)$, the relation between time orientation, future pointed observer, time inversion operation and past pointed observers is not the usual one.

\subsection{\small{Examples of Finsler spacetimes}}
We collect several examples of Finsler spacetimes investigated in the literature.
 The examples below do not exhaust the intense use of Finsler
geometries in physical applications. On the other hand, such a bunch of examples partially
motivates the mathematical investigation of Finsler spacetimes.
\begin{ejemplo}
The first example to consider are Lorentzian spacetimes $(M,h)$, where $h$ is a Lorentzian metric.
 In this case, the Lagrangian is given by
 \begin{align}
 L(x,y)=\,h_x(y,y),\quad y\in T_xM.
 \end{align}
\end{ejemplo}
\begin{ejemplo} Let $M$ be an $n$-dimensional manifold and let us consider the following Lagrangian function,
\begin{align}
L(x,y)=\frac{1}{2}\Big(\ell(x,y)^2\,-U_i(x)U_j(x)\,y^i\,y^j\Big)
\end{align}
\label{ejemplofinsler1}
where $U(s)$ defines a $1$-form on M  and $\ell(x,v)$ are such that
following conditions:
\begin{enumerate}
\item $ \ell(x,ky)=\,k\,\ell(x,y)$ for positive $k$,
\item $\frac{\partial ^2\ell^2 (x,y)}{\partial y^i\partial y^j}
\,w^i\,w^j\,>0$ if $U(x)(w)>0$ and
\item There is a unique vector field $V(x)$ defined by $U(V)=-1$ and
$\frac{\partial^2 \ell^2 (x,y)}{\partial y^i\partial y^j}V^i(x)=0$, $\omega\neq 0$.
\end{enumerate}
These conditions guarantee that the matrix of fundamental tensor components
\begin{align}
g_{ij}(x,y)=\,\frac{1}{2}\,\frac{\partial^2\ell^2(x,y)}{\partial y^i\partial y^j}-\,U^iU^j
\end{align}
is non-degenerate and with signature $(-1,1,...,1)$.
The relevance of this example to physics
resides in that it describes light propagation in a linear, dielectric and permeable medium \cite{Perlick2000}.
\end{ejemplo}
\begin{ejemplo}
A family of Finsler spacetimes that have been considered in the physics literature
 are based on Berwald-Moor Finsler metrics \cite{LämmerzahlLorekDittus}.
Let $(M,\eta)$ be the Minkowski spacetime and $W$ a timelike vector field on $M$. An Euclidean metric induced by $W$
is
\begin{align}
\hat{\eta}_x(y,y)=\eta_x(y,y)-\,2\frac{\eta^2_x(y,W)}{\eta_x(W,W)}.
\end{align}
Let $\hat{y}$ be the orthogonal component of $y$ to $W$ using $\hat{\eta}$ and $\phi$ a $2p$-tensor.
Then the fundamental tensor $g$ is of the form
\begin{align}
g_x(y,y)=\,\eta_x(y,y)+\hat{\eta}_x(y,y)\,\Big(\frac{1}{p}\frac{\phi(\hat{y},....,\hat{y})}{\hat{\eta}_x(y,y)^p}\Big).
\label{ejemplolammerzahl}
\end{align}
This tensor determines a Finsler spacetime iff $\phi$ is small enough compared with $\hat{\eta}$.
Experimentally, Finsler spacetimes of Berwald-Moor type are constrained to be Lorentzian with a very high accuracy \cite{LämmerzahlLorekDittus}.
\label{ejemploLLD}
\end{ejemplo}
\begin{ejemplo}
Not directly related with physical models is the following example \cite{Beem1}. The spacetime manifold is $M=R^3$
and the Lagrangian is the highly  non-reversible function $L$
\begin{align}
L(y)=\,\frac{(y^1)^3\,-y^1(y^2)^2}{\big((y^1)^2+\,(y^2)^2\big)^{\frac{1}{2}}}.
\end{align}
Then $L(-y)=\,-L(y)$ and the indicatrix has six connected components.
\label{ejemplobeem}
\end{ejemplo}

The following examples share the common fact that they are not regular
 in the whole slit tangent space of a spacetime manifold. In order to consider such examples one needs to relax the conditions of the Finsler spacetime
 (some notions of {\it weak Finsler structures} can be found in \cite{Randers, PfeiferWohlfarth}).
 \begin{ejemplo}
 Rutz has investigated a non-Riemannian solutions of a Einstein-Finsler theory in vacuum \cite{Rutz}.
Let coordinates $(t,r,\theta, \varphi)$ be local spherical coordinate system.  In spherical coordinates, a tangent vector $y\in\,TM$ is expressed as
 \begin{align*}
y=\,y_t\frac{\partial}{\partial t}+\,y_r\frac{\partial}{\partial r}+\,y_{\theta}\frac{\partial}{\partial t}+
\, y_{\varphi}\frac{\partial}{\partial\, \varphi},
 \end{align*}
Ruth's Finsler spacetime is 
a static, spherical symmetric, Finsler space-time, with spacetime manifold $M=\,R\times R^+\times S^2$ and Lagrangian
\begin{align}
L(x,y):=\Big(-\Big(1-\frac{2m}{r}\Big)
\Big(1-\,\delta\,\frac{d\Omega}{dt}\Big)\,dt^2+
\frac{1}{\Big(1-\frac{2m}{r}\Big)\,}dr^2\,+r^2d\Omega^2\Big)\cdot(y,y), 
\label{rutzexample}
\end{align}
where the function $\frac{d\Omega}{dt}$ is defined by
 \begin{align*}
\frac{d\Omega}{dt}:=\,\frac{\sqrt{y^2_{\theta}+\,\sin^2\theta\,y^2_{\varphi}}}{y_t},\,\,y_t\neq 0
\end{align*}
and the parameter $\delta$ is small compared with $1$.

The Schwarzschild's
solution of Einstein's equations \cite{BEE},
\begin{align}
L_S(x,y):=\,ds^2_{S}\cdot(y,y)=\,\Big(-\Big(1-\frac{2m}{r}\Big)\,dt^2+
\frac{1}{\Big(1-\frac{2m}{r}\Big)\,}dr^2\,+r^2d\Omega^2\Big)\cdot(y,y),
\label{Schwarzschild}
\end{align}
 $(M,L_R)$ is a {\it singular} Finsler spacetime, since it is not regular in the full $N$. Indeed,  it is smooth on $N\setminus A$, where $A$ contains the set
 where $\frac{d\Omega}{dt}=0$, the sub-manifolds $\{y\in T_xM\,|\,y_t=0\}\hookrightarrow T_xM$ and at the Schwarzschild radius $r_S\,=2m$.
 The corresponding fundamental tensor
 \begin{align*}
 (g_R)_{ij}=\, \frac{\partial^2 L_R(x,y)}{\partial y^i \partial y^j}
 \end{align*}
  is non-degenerate on $N\setminus A$ and has Lorentzian signature for $\delta$ small enough\footnote{Note that since the fundamental tensor associated with $L_R$ is different than the {\it generalized metric} $g_R$ associated with the line element $ds^2_R$. That is, the form \eqref{rutzexample} defines a Finsler spacetime in the sense of Beem's and a generalized metric in the sense of Miron and Anastasiei \cite{MironAnastasiei}. The distance function coincide, but the {\it fundamental tensors} are different.}. Also, note that is not a reversible metric.
\end{ejemplo}
\begin{ejemplo}
The {\it rainbow metric} is a phenomenological description of the modification of the dispersion relations
produced by possible quantum gravity corrections \cite{GirelliLiberatiSindoni}.
Let $(M,\eta)$ be a stationary Lorentzian spacetime such that $\mathcal{L}_W \eta=0$.
There is a foliation on $M$ given by integral curves of $W$.
The orthogonal spacelike hypersurfaces $\Sigma_t$ furnish an induced Riemannian metric $\bar{\eta}$ by isometric embedding.
The rainbow metric is determined by the following lagrangian function $L$ (compare with \cite{GirelliLiberatiSindoni}),
\begin{align}
L(x,y)=\,\Big(\sqrt{\eta(x,y)}-\,C_1(m)\frac{\bar{\eta}(\bar{y},\bar{y})^{\frac{3}{2}}}{\eta^{\frac{1}{2}}(x,y)}\Big)^2.
\label{rainbowmetric}
\end{align}
This metric is not regular in the light cone $\eta(x,y)=0$. This singularity is related with the mass of the particle $m$.
Therefore, each specie of elementary particle has its {\it particular metric}
(this is why the name rainbow metric). The rainbow metric is non-reversible.
\end{ejemplo}
\begin{ejemplo}
Related with very special relativity of Cohen and Glashow, there are related Finsler metrics of Bogoslovsky type.
In particular, very special relativity group leaves invariant the line element of norm
\begin{align}
L(x,y)=\,\eta(y,y)^{(1-b)}\,n (y)^{2b},
\label{bogosmetric}
\end{align}
with $\eta$ the Minkowski metric in $4$ dimensions, $n$ a $1$-form
corresponding to the null direction $\nu^{\mu}=\delta^\mu_+$ and $b\neq 1$ is
the deformation parameter \cite{GibbonsGomisPope}.
Bogoslovsky metric contains singularities on the cone $\eta(y,y)=0$.
\end{ejemplo}
\begin{ejemplo}
Bi-metric theories have been considered in the literature associated with birefringent crystal optics \cite{SkakalaVisser}.
They are constructed from two Lorentzian spacetimes  $(M,L_+)$ and $(M,L_-)$ by the Lagrangian
\begin{align}
L(x,y):=\,\sqrt{L_+(x,y)\,L_-(x,y)}.
\end{align}
These metrics are singular on each of the null cones $L_+(x,y)=0$ and $L_-(x,y)=0$,
since the corresponding fundamental tensor $g_{ij}=\frac{1}{2}\,\partial_{y^i}\partial_{y^j}\,L(x,y)$ is not smooth \cite{SkakalaVisser2}.
\end{ejemplo}

\subsection{\small{Variational setting}}
\begin{definicion}
An affine parameterized geodesic of a Lagrangian $L$ is a solution of the Euler-Lagrange equation
\begin{equation}
\frac{d}{ds}\,\frac{\partial L(\lambda(s),\dot{\lambda}(s))}{\partial \dot{\lambda}^i}-\frac{\partial
L(\lambda(s),\dot{\lambda}(s))}{\partial \lambda^i}=0,\quad i=1,...,n,
\label{affineparameterizedgeodesic}
\end{equation}
with $\dot{\lambda}^i(s)=\,\frac{d\lambda^i(s)}{ds}$. In this case $s$ is an affine parameter.

An arbitrarily parameterized geodesic is a solution of the differential equation
\begin{equation}
\frac{d}{ds}\,\frac{\partial L(\lambda(s),\dot{\lambda}(s))}{\partial \dot{\lambda}^i}-\frac{\partial
L(\lambda(s),\dot{\lambda}(s))}{\partial \lambda^i}=\,f(s)\dot{\lambda}^i,\quad i=1,...,n
\label{affineparameterizedgeodesic2}
\end{equation}
for a given function $f:I\longrightarrow M$ and with $\dot{\lambda}^i(s)=\,\frac{d\lambda^i(s)}{ds}$.
\end{definicion}
Given an arbitrarily parameterized geodesic of an  affine connection on $M$,
 it is possible to find a positive re-parameterization such that with the new parameter the curve is an affine geodesic.

 Using the equation (\ref{affineparameterizedgeodesic}) and the homogeneity condition (\ref{homogeneouscondition}),
 one can show that $L(\lambda,\dot{\lambda})$ is preserved along affine parameterized geodesics,
\begin{align*}
\frac{d}{ds}L(\lambda(s),\dot{\lambda}(s)) & =\,\frac{\partial L(\lambda(s),\dot{\lambda}(s))}{\partial \lambda^i}\,
\dot{\lambda^i}\,+\frac{\partial L(\lambda(s),\dot{\lambda}(s))}{\partial \dot{\lambda}^i}\,\ddot{\lambda^i}\\
& =\,\frac{d}{ds}\Big(\frac{\partial L(\lambda(s),\dot{\lambda}(s))}{\partial \dot{\lambda}^i}\Big)\dot{\lambda}^i\,
+\frac{\partial L(\lambda(s),\dot{\lambda}(s))}{\partial \dot{\lambda}^i}\,\ddot{\lambda}^i\\
& =\,\frac{d}{ds}\Big(\frac{\partial L(\lambda(s),\dot{\lambda}(s))}{\partial \dot{\lambda}^i}\,\dot{\lambda}^i\Big)\\
& =\,2\frac{d}{ds}\Big(L(\lambda(s),\dot{\lambda}(s))\Big),
\end{align*}
from which follows that $L(\lambda,\dot{\lambda})$ is constant along $\lambda$
 (and therefore, also along any equivalent arbitrarily re-parameterized geodesic).
Therefore, a causal geodesic is a geodesic with $g_{\dot{\lambda}(s)}(\dot{\lambda},\dot{\lambda})\leq 0$;
 for a timelike geodesic $g_{\dot{\lambda}(s)}(\dot{\lambda},\dot{\lambda})\leq 0$ and for a
lightlike geodesic $g_{\dot{\lambda}(s)}(\dot{\lambda},\dot{\lambda})=0$. Note that the causal character
of a geodesic is preserved by re-parameterization and that time orientation is preserved by monotone increasing re-parameterizations.

Let us consider a point $q\in M$, a constant $c\leq 0$ and a future pointed, timelike curve $\gamma:I\longrightarrow M$.
 Then the {\it space of admissible curves} is the space
\begin{align*}
\mathcal{C}_{q,\gamma,c}:=\Big\{\,& \lambda:[0,1]\longrightarrow M, \quad \textrm{smooth such that}\\
& 1. \,\lambda(0)=q,\\
& 2. \,\exists \,\tau(\lambda)\in I\,\,s.t.\,\,\lambda(1)=\,\gamma(\tau(\lambda)),\\
& 3. \,L(\lambda(s),\dot{\lambda}(s))=-c^2, \,\forall s\in [0,1],\\
& 4. \, g_{ij}(\lambda(s),\dot{\lambda}(s))\,\dot{\lambda}^i(s)\,T^j(\lambda(s))<0\,\Big\}.
\end{align*}
Note that if $\lambda\in \mathcal{C}_{q,\gamma,c}$, $\lambda$ will not be parameterized necessarily by the {\it proper time}, defined by the integral
\begin{align}
t_{\tilde{\lambda}}(\tilde{\lambda}(r))=\int^{r_2}_{r_1}\,\sqrt{-g_{ij}(\tilde{\lambda}(r),
\dot{\tilde{\lambda}}(r))\,\dot{\tilde{\lambda}}^i(r)\,\dot{\tilde{\lambda}}^j(r)}\,dr.
\label{propertime}
\end{align}
\begin{definicion}
An allowed variation of $\lambda\in\,\mathcal{C}_{q,\gamma,c}$ is a smooth map
\begin{align*}
\Lambda :(-\epsilon_0,\epsilon_0)& \times\, [0,1] \longrightarrow M,\,\,\,\epsilon_0 >0
\end{align*}
such that
\begin{enumerate}
\item Each of the curves $\Lambda(\epsilon,\cdot)$ is allowed,
\begin{align*}
(\epsilon, s)\mapsto \Lambda(\epsilon,\cdot)\in \,\mathcal{C}_{q,\gamma,c}\quad \forall \epsilon \in [-\epsilon_0,\epsilon_0],
\end{align*}
\item The central curve is $\Lambda(0,s)=\lambda(s)$.
\end{enumerate}
\label{alowedvariation}
\end{definicion}
We introduce two functionals relevant for our purposes,
\begin{definicion}
Let $\mathcal{C}^{\infty}([0,1],M)$ be the space of smooth parameterized curves of
$M$ parameterized in the interval $[0,1]$. The energy functional is
\begin{align}
E: \mathcal{C}^{\infty}([0,1],M)\longrightarrow R,\quad \lambda \mapsto E(\lambda):=\int^1_0\,L(\lambda(s),\dot{\lambda}(s))\,ds,
\label{energyfunctional}
\end{align}
\end{definicion}
Note that for any allowed variation, the energy $E$ of each curve is $-c^2$. Therefore, we
are considering causal curves with prescribed energy. As a consequence of the prescription of the energy one has that
\begin{align}
\frac{d}{d\epsilon}\Big|_{\epsilon=0}\,\Big(\int^1_0\,L(\Lambda(\epsilon,s),\dot{\Lambda}(\epsilon, s))\,ds\Big)=0.
\label{criticalconditioforE}
\end{align}
All the curves in the class $\mathcal{C}_{q,\gamma, c}$ have constant energy equal to $E=-c^2$.
\begin{definicion} Let $\mathcal{C}_{q,\gamma,c}$ be the space of admissible curves.
The time arrival functional is
\begin{align}
\tau:\mathcal{C}_{q,\gamma,c}\longrightarrow R,\quad
& \lambda \mapsto \tau(\lambda).
\label{timearrivalfunctional}
\end{align}
\end{definicion}
\section{\small{Fermat's principle for causal curves in time oriented Finsler spacetimes}}
\subsection{\small{Regularity of the time arrival functional}}
In standard treatments of Fermat's principle for lightlike geodesics it is assumed that the
time arrival functional acting on any allowed variation $\Lambda(\epsilon,s)$ is of class
$\mathcal{C}^1$ in the variable $\epsilon$ \cite{Perlick, Perlick06}.
Such regularity  holds when $L$ is a Lorenztian metric
and the allowed curves $\lambda$ are timelike \cite{GiannoniMasielloPiccione}. Indeed one has the following result,
\begin{proposicion}
Let $(M,L)$ be a Finsler spacetime, $\Lambda: (-\epsilon,\epsilon)\times \, [0,1]\longrightarrow M$ a variation of a causal geodesic $\lambda$ and
$\gamma:I \longrightarrow M$ a timelike, positive temporary oriented curve.
Then the function $\tau(\Lambda(\epsilon,\cdot))$ is smooth on $\epsilon$.
\label{proposicionon regularity of tau}
\end{proposicion}
{\bf Proof}. Let us consider the time arrival functional acting on the variation
$\Lambda(\epsilon,s)$, i.e., the function
\begin{align*}
t:(-\epsilon_0,\epsilon_0 )\longrightarrow R,\quad\epsilon\mapsto \gamma^{-1}(\Lambda(\epsilon,1))=\gamma^{-1}\circ \Lambda(\epsilon,1),
\end{align*}
The function
\begin{align*}
\Lambda(\cdot,1):(-\epsilon_0,\epsilon_0 )\longrightarrow M,\quad\epsilon\mapsto \Lambda(\epsilon,1)
\end{align*}
is smooth. Since $\gamma(\sigma)$ is smooth, $(\gamma^{\mu}) '(\sigma)=\frac{d\gamma^{\mu} (\sigma)}{d\sigma}\neq 0$
 for any $\sigma\in\,I$ and $\gamma(\sigma)$
does not have self-intersections, $\gamma^{-1}:\gamma(I)\longrightarrow R$ is smooth.
Therefore, since $t(\epsilon)=\,\gamma^{-1}(\Lambda(\epsilon,\cdot))=\,\tau(\Lambda(\epsilon,\cdot))$
 is smooth on $\epsilon$ the result follows. \hfill$\Box$

The smoothness on $\epsilon$ of the time arrival functional is fundamental in the formulation of Fermat's principle as well as for related results.

\subsection{\small{Fermat's principle in Finsler spacetimes for causal curves}}
Let us fix the time positive oriented timelike curve $\gamma:I\longrightarrow M $.
Fermat's principle for causal curves can be stated as follows
\begin{proposicion}
Let $(M,L)$ be a time orientable Finsler spacetime. Then the causal curve $\lambda:[0,1]\longrightarrow M$ is
a geodesic (pre-geodesic in the lightlike case) of $L$ iff it is a critical point of the time arrival functional (\ref{timearrivalfunctional}),
\begin{align}
\frac{d}{d\epsilon}\Big|_{\epsilon=0}\tau(\Lambda(\epsilon,s))=0,
\end{align}
for any allowed variation $\Lambda(\epsilon,s)$ of $\lambda(s)$.
\label{fermatprinciplefortimelikegeodesics}
\end{proposicion}
\begin{comentario}
This is a generalization of the Finslerian version of Fermat's principle for lightlike curves  obtained in \cite{Perlick} and of the Lorentzian
 Fermat's principle for timelike curves
\cite{GiannoniMasielloPiccione}. Note that the allowed curves are different from the above mentioned principles.
For instance, one does not require time orientation for the light-like curves \cite{Perlick};
for Fermat's principle contained in \cite{GiannoniMasielloPiccione} the notion of time-orientation
is slightly different than the principle considered in this work. Also, because a technicality in the prove of the theorem, we will require future pointed oriented curves. This is in contrast with \cite{Perlick}, where only a positivity orientation is required when observed by $\gamma$.
\end{comentario}
\begin{comentario}
For a timelike geodesic, the parameter of $\lambda(s)$
 is an affine parameter. This is not the case if $ \lambda$ is a lightlike curve.
\end{comentario}
Before we prove {\it proposition} \ref{fermatprinciplefortimelikegeodesics},
 let us write some intermediate formulas. First note that a smooth curve $\lambda$
is a critical point of the functional energy $E$ iff equation \eqref{criticalconditioforE} holds
for any allowed variation $\Lambda(\epsilon,s)$. Also note that since $\Lambda(\epsilon,1)=\,
\gamma(\tau(\Lambda(\epsilon, 1)))$ one has that in local coordinates
\begin{align}
\frac{d}{d\epsilon}\Big|_{\epsilon=0}\,\Lambda^i(\epsilon,1)=\,({\gamma^i})'(\tau(\lambda))\,
 \frac{\partial}{\partial\epsilon}\Big|_{\epsilon=0}\,\big(\tau(\Lambda(\epsilon,s))\big),\quad i=1,...,n,
\label{relationLambdatau}
\end{align}
for any allowed variation $\Lambda(\epsilon,s)$.

{\bf Proof of proposition \ref{fermatprinciplefortimelikegeodesics}}. The ``only if'' is proven following a similar argument
as in \cite{Perlick06}.  The condition that all the curves in the allowed variation $\Lambda$ are of fixed energy and that
the allowed variation $\Lambda$ is indeed a smooth function on $\epsilon$ and $s$ implies the following relation,
\begin{align*}
0 & =\,\int^1_0\,\frac{\partial}{\partial\epsilon}\big|_{\epsilon=0}\,\big(L(\Lambda(\epsilon,s),\dot{\Lambda}(\epsilon, s))\big)\,ds\\
& =\,\int^1_0\,\Big(\frac{\partial L}{\partial \lambda^i}\,\frac{\partial}{\partial\epsilon}\big|_{\epsilon=0}(\Lambda^i(\epsilon,s))\,
+\frac{\partial L}{\partial \dot{\lambda}^i}\,\frac{\partial}{\partial\epsilon}\big|_{\epsilon=0}\dot{\Lambda}^i(\epsilon,s)\Big)\,ds\\
& =\,\int^1_0\,\Big(\frac{\partial L}{\partial \lambda^i}\,\frac{\partial}{\partial\epsilon}\big|_{\epsilon=0}(\Lambda^i(\epsilon,s))\,
+\frac{\partial L}{\partial \dot{\lambda}^i}\,\frac{d}{ds}\frac{\partial}{\partial\epsilon}\big|_{\epsilon=0}{\Lambda}^i(\epsilon,s)\Big)\,ds\\
& = \,\int^1_0\,\Big(\Big(\frac{\partial L}{\partial \lambda^i}\,-\frac{d}{ds}\,\frac{\partial L}{\partial \dot{\lambda}^i}\Big)
\frac{\partial}{\partial\epsilon}\big|_{\epsilon=0}\Lambda^i(\epsilon,s))\Big)\,ds\,+\,
\Big(\frac{\partial L(\lambda,\dot{\lambda}(s))}{\partial \dot{\lambda}^i}\,\frac{\partial}{\partial\epsilon}\big|_{\epsilon=0}\Lambda^i\Big)\Big|^1_0.
\end{align*}
Then using the relation (\ref{relationLambdatau}) one obtains
\begin{multline*}
\int^1_0\,\Big(\frac{\partial L}{\partial \lambda^i}\,-\frac{d}{ds}\,\frac{\partial L}{\partial \dot{\lambda}^i}\Big)
\frac{\partial}{\partial\epsilon}\big|_{\epsilon=0}(\Lambda^i(\epsilon,s))\,ds+ \frac{\partial L(\lambda,\dot{\lambda})}{\partial \dot{\lambda}^i}\,({\gamma^i})'(\tau(\lambda))\frac{d}{d\epsilon}\big|_{\epsilon=0}\,
\big(\tau(\epsilon,s)\big)=0
\end{multline*}
and by the homogeneity property of $L$,
\begin{multline}
\int^1_0\,\Big(\frac{\partial L}{\partial \lambda^i}\,-\frac{d}{ds}\,\frac{\partial L}{\partial \dot{\lambda}^i}\Big)
\frac{d}{d\epsilon}\big|_{\epsilon=0}(\Lambda^i)\,ds+\big(g_{jl}(\lambda(1),{\lambda}'(1))\dot{\lambda}^l(1)({\gamma}^j)'(\tau(\lambda))\big)\,\frac{d}{d\epsilon}
\big|_{\epsilon=0}\,\big(\tau(\epsilon,s))
=0.
\label{relationeulerlagrangetau}
\end{multline}
Given a curve $\lambda \in \mathcal{C}_{q,\gamma,c}$, one has the condition
\begin{align}
g_{\dot{\lambda}(s)}(\dot{\lambda}(1),\gamma'(\tau(\lambda)))=\,g_{jl}(\lambda(1),\dot{\lambda}(1))\dot{\lambda}^l(1)({\gamma}^j)'(\tau(\lambda))\neq 0.
\label{causalitycondition}
\end{align}
Let us parameterize the geodesic by an affine parameter $s\in[0,1]$, which means that equation (\ref{affineparameterizedgeodesic}) holds.
 Then it is clear from (\ref{relationeulerlagrangetau}) that for curves of fixed energy $E=-c^2\,\leq\,0$,
 the solutions of the Euler-Lagrange equations are critical points of the time arrival functional,
\begin{align*}
\frac{\partial L(\lambda,\dot{\lambda})}{\partial \lambda^i}\,-\frac{d}{ds}\,\frac{\partial L(\lambda,
\dot{\lambda})}{\partial \dot{\lambda}^i}=0,\quad E=-c^2\,\Rightarrow\frac{d}{d \epsilon}\big|_{\epsilon=0}\tau(\Lambda(\epsilon,s))=0.
\end{align*}
This implication is independent of the signature of the metric. It strongly depends on the requirement
that the energy $E$ has a fixed valued for all the allowed curves.

The ``if" implication in proposition \ref{fermatprinciplefortimelikegeodesics} can be proved as follows.
Let us consider a variation of $\mathcal{C}_{q,\gamma,c}\ni\,\lambda:[0,1]\longrightarrow M$ defined by
\begin{align}
\Lambda^{\alpha}(\epsilon,s)=\lambda^{\alpha}(s)\,+\epsilon A^{\alpha}(s),\quad \alpha=1,...,n-1,
\label{definitionA}
\end{align}
with $A^{\alpha}$ arbitrary smooth functions. In order to be an allowed variation along the curve
$\lambda(s)\in\,\mathcal{C}_{q\gamma,c}$, $\Lambda(\epsilon,s)$ must satisfy  the condition
\begin{align}
L(\Lambda,\dot{\Lambda})=-c^2
\label{equationforLambdaN}
\end{align}
with initial condition $\Lambda(\epsilon,0)=-c^2$. If we cover the image of $\Lambda$ by local coordinate charts, this condition can be expressed explicitly in a convenient local coordinate system such that the last coordinate
$x^n$ corresponds to the integral curves of the time orientation vector field $T(x)$. In such coordinate system
$\sum^n_{i=1}\,g_{ni}(\Lambda,\dot{\Lambda})\dot{\Lambda}^i\,\neq 0$ holds and one has the constrain
\begin{align}
c^2+\,\sum^{n-1}_{\alpha,\beta =1}\,g_{\alpha\beta}(\Lambda,\dot{\Lambda})\dot{\Lambda}^{\alpha}\dot{\Lambda}^{\beta}+
\,\sum^n_{i=1}\,g_{ni}(\Lambda,\dot{\Lambda})\dot{\Lambda}^i\,\dot{\Lambda}^n\,=0.
 \label{rewritingequationforLambdaN}
\end{align}
The fact that $\sum^n_{i=1}\,g_{ni}(\Lambda,\dot{\Lambda})\dot{\Lambda}^i\,\neq 0$ can be used to solve $\dot{\Lambda}^n$
  in \eqref{rewritingequationforLambdaN},
\begin{align}
\frac{1}{
\,\sum^n_{i=1}\,g_{ni}(\Lambda,\dot{\Lambda})\dot{\Lambda}^i}\,c^2+\,\frac{\sum^{n-1}_{\alpha,\beta =1}
\,g_{\alpha\beta}(\Lambda,\dot{\Lambda})\dot{\Lambda}^{\alpha}\dot{\Lambda}^{\beta}}{
\,\sum^n_{i=1}\,g_{ni}(\Lambda,\dot{\Lambda})\dot{\Lambda}^i}+\,\dot{\Lambda}^n\,=0.
\label{rewritingequationforLambdaN2}
\end{align}
Let us fix the value of $\epsilon$. Since \eqref{rewritingequationforLambdaN2} does not contain any derivative $\frac{\partial}{\partial\epsilon}$, the variable
 $\epsilon$ can be considered as a continuous parameter of the differential
 equation (\ref{rewritingequationforLambdaN2}) that s considered as an ODE. Then one
uses standard ODE theory to establish local existence, uniqueness and smoothness of $\Lambda(\epsilon,s)$ on the parameter $\epsilon$ for
$s\in [0,s_0]$ for some $s_0$ (see for instance \cite[Chapter~1]{Chicone}). Such solution can be extended further. Indeed, the solution
can be extended to an interval $]s_0-\delta,s_0+\delta[$, with the convenient initial
data at the point $s_0-\delta\in \,[0,s_0[$ and with a smooth dependence on $0<\epsilon <\epsilon_1\leq \epsilon_0$. Repeating this procedure
one can extend the solution to a finite collection of open sets of $R$ which
is maximal and contained in $\lambda([0,1])$. Let $s_{max}\leq \,1$ be the maximal value of $s$
such that the dependence on $\epsilon \in\, \Lambda(\epsilon,s)$ is smooth. There are two possibilities,
\begin{enumerate}
\item $\lambda\in \, \mathcal{C}_{q,\gamma,c}$. In this case, $s_{max}=1$, since otherwise one
can extend $s_{max}$ leading to a contradiction with ODE theory.
\item  $\lambda$ does not intersect $\gamma$. This is in contradiction with the hypothesis $\lambda(1)=\gamma(\tau(\lambda))$.
\end{enumerate}
Let us consider the minimum $\epsilon_{min}$ of the above $\epsilon_k$, that by compactness of $[0,1]$ must be finite.
Using again a compactness argument, one can show that there is a finite open cover of $\lambda:[0,1]\longrightarrow M$
such that the differential equation (\ref{rewritingequationforLambdaN2}) contains an unique solution on
each local chart. Because paracompact property of $M$, there is an adapted partition of the unity.
Using bump functions \cite{Warner}, one can patch in a smooth way a solution in an open subset of $[0,1]$,
obtaining a global solution $\lambda$ for (\ref{rewritingequationforLambdaN2}) in $[0,1]$. This proves the existence of allowed
variations of $\lambda$ satisfying the ODE (\ref{rewritingequationforLambdaN2}).

Using the local existence and uniqueness of ODE's theory locally one can write the expression
\begin{align}
\dot{\Lambda}^n(\epsilon,s)(\epsilon,s)=\,\tilde{f}_0(s,A(s))+\epsilon \,\int^s_0\,\big(\,\sum^{n-1}_{\alpha=1}\,A^{\alpha}(s)\tilde{h}_{\alpha}(s)
+\,\epsilon^2\,\tilde{f}_2(s,A(s),\bar{\epsilon})\big) \,ds
\label{relationdotLambdaepsilon}
\end{align}
for some unique, smooth functions $\{\tilde{h}_{\alpha}(s)\}$.
 Integrating respect to $s$ both sides (\ref{relationdotLambdaepsilon}) one obtains
\begin{align*}
{\Lambda}^n(\epsilon,s)(\epsilon,s)=\,{f}_0(s,A(s))+\epsilon \,\big(\,\sum^{n-1}_{\alpha=1}\,A^{\alpha}(s)\tilde{h}_{\alpha}(s)\big)
+\,\epsilon^2\,{f}_2(s,A(s),\bar{\epsilon}),
\end{align*}
Then
the variation of the relation (\ref{relationeulerlagrangetau}) is equivalent to
 \begin{align*}
\int^1_0\,\Big(\sum^{n-1}_{\alpha=1}\,\big(\frac{\partial L}{\partial \Lambda^{\alpha}}-&\,
\frac{d}{ds}\frac{\partial L}{\partial \dot{\Lambda}^{\alpha}}\big)\,{A}^{\alpha}(s) +\,
\big(\frac{\partial L}{\partial {\Lambda}^n}-\,\frac{d}{ds}\frac{\partial L}{\partial \dot{\Lambda}^n}\big)\,
\big(\,\sum^{n-1}_{\alpha=1}\,A^{\alpha}(s)\tilde{h}_{\alpha}(s)\big)\Big)\,ds\\
& = -\int^1_0\,\frac{d}{ds}\Big(\,({\gamma}^i)'(\tau(\lambda))\,\frac{\partial L(\lambda,\dot{\lambda}(s))}{\partial
\dot{\lambda}^i}\,
\frac{d}{d\epsilon}|_{\epsilon=0}\big(\tau(\Lambda(\epsilon,s))\Big)\,ds.
 \end{align*}
 For critical points of the arrival time functional it holds that
 \begin{align*}
 \frac{d}{d\epsilon}\big|_{\epsilon=0}(\tau(\Lambda(\epsilon,\cdot)))=0.
 \end{align*}
 This implies
  \begin{multline*}
\int^1_0\,\Big(\sum^{n-1}_{\alpha=1}\,\big(\frac{\partial L}{\partial \Lambda^{\alpha}}-\frac{d}{ds}
\frac{\partial L}{\partial \dot{\Lambda}^{\alpha}}\big)\,{A}^{\alpha}(s) +\,\big(\frac{\partial L}{\partial {\Lambda}^n}-\,
\frac{d}{ds}\frac{\partial L}{\partial \dot{\Lambda}^n}\big)\,\big(\,\sum^{n-1}_{\alpha=1}\,A^{\alpha}(s)\tilde{h}_{\alpha}(s)\big)\Big)\,ds=0.
 \end{multline*}
Therefore, one can write the relations
 \begin{align}
0=\,\int^1_0\,\Big(\sum^{n-1}_{i=1}\frac{\partial L}{\partial \lambda^i}\,-\frac{d}{ds}\,\frac{\partial L}{\partial \dot{\lambda}^i}\Big)
\,B^i(s)\Big)\,ds
\label{newrelationeulerlagrangetau}
\end{align}
 for the arbitrary, small enough functions $(B^1(s),...,B^{n-1}(s))$.
 Also note that since for each $s$ for some $(A^1(s),...,A^{n-1}(s))$ the expression
 \begin{align*}
 \big(\,\sum^{n-1}_{\alpha=1}\,A^{\alpha}(s)\tilde{h}_{\alpha}(s)\big)\Big),
 \end{align*}
 is arbitrarily large, one finds the condition
 \begin{align}
 \frac{\partial L}{\partial {\lambda}^n}-\,
\frac{d}{ds}\frac{\partial L}{\partial \dot{\lambda}^n}=0.
\label{lagrangeequationforn}
 \end{align}
 The functions $\{B^1(s),...,B^{n-1}(s)\}$ determine a $(n-1)$-dimensional
 vector space at each point $\lambda(s)$. Each of such $(n-1)$-vector subspaces
 are orthogonal to the tangent vector $\dot{\lambda}(s)$. Let us consider
 a basis  $\{B_1(s),...,B_{n-1}(s)\}$ for the subspace of $T_{\lambda(s)}M$ orthogonal to $\dot{\lambda}$.
Using the fact that $g$ is non-degenerate, the expression (\ref{newrelationeulerlagrangetau}) can be rewritten as
\begin{align}
0=\,\int^1_0\,g_{\dot{\lambda}}(\omega(s),B_{\alpha}(s))\, ds
\label{omega}
\end{align}
for some vector field $\omega:[0,1]\longrightarrow M$ defined by the relation
\begin{align*}
g_{\dot{\lambda}}(\omega(s),B(s))=\big(\frac{\partial L}{\partial \Lambda^{i}}-\frac{d}{ds}
\frac{\partial L}{\partial \dot{\Lambda}^{i}}\big)\,{B}^{i}(s)
\end{align*}
for each space-like vector $(B_1(s),...,B_{n-1}(s))$.
Since the functions $B^i(s)$ are arbitrary, we can extract from the integral the
local condition
\begin{align}
g_{\dot{\lambda}}(\omega(s),B_{\alpha}(s))=0.
\label{relationbetweenomegaB}
\end{align}
This implies that for a fixed $s$ the vector $\omega(s)$ is parallel to $\dot{\lambda}(s)$. This can be proven as follows.
If $\omega(s)$ is not parallel to $\dot{\lambda}(s)$,
there is vector $C(s)$ such that $g_{\dot{\lambda}}(C(s),B(s))=0$ for all $B(s)\in\, span\{B_1,...,B_{n-1}\}$.
 As we say, the dimension of the variational fields $B(s)$ orthogonal to $\dot{\lambda}(s)$ is $n-1$ for all $s$.
The dimension of $span\{\dot{\lambda}(s),B_i(s),C(s)\}$ must be maximum $n$. Therefore, $C(s)=0$ and it follows $\omega(s)$ is parallel to $\dot{\lambda}(s)$.
Since equation \eqref{relationbetweenomegaB} holds for the arbitrary $n-1$ $B(s)$ functions, it must hold that $\omega^i(s)=0,\,i=1,...,n-1$.
\hfill$\Box$

It is remarkable the use of time orientation in the above proof to isolate $\dot{\Lambda}(\epsilon,s)$ and be able to use ODE theory.  Also note, that as in the positive case, if $\gamma:I\to M$ with $I=[a,b]$ compact is a geodesic, then $\hat{\lambda}:-I\to M, \,\hat{\gamma}:=\gamma(-s)$ is not necessarily a geodesic. Therefore, if $\gamma$ is a critical point of the time arrival functional $\tau$, then the reversed curve $\hat{\gamma}$ is not a geodesic.
\subsection{\small{Formula for the first variation of the time arrival functional}}
The above calculation provides the formula of the first variation of the time arrival functional. Writing
\begin{align*}
g_{jl}(\Lambda(1,\epsilon),\dot{\Lambda}(1,\epsilon)\dot{\Lambda}^l(1,\epsilon)({\gamma}^j)'(\tau(\Lambda(\epsilon,s)))
=g_{\dot{\Lambda}}(\dot{\Lambda},{\gamma}'(\tau(\Lambda))),
\end{align*}
 the first variation of the functional $\tau(\Lambda)$ is given by the expression:
\begin{multline}
g_{\dot{\Lambda}}(\dot{\Lambda},{\gamma}'(\tau(\Lambda)))\,\frac{d}{d\epsilon}\,\big(\tau(\Lambda(\epsilon,s))\big)
\\=-\int^1_0\,\Big(\frac{\partial L(\Lambda(\epsilon,s),\dot{\Lambda}(\epsilon,s))}{\partial \Lambda^i}\,
-\frac{d}{ds}\,\frac{\partial L(\Lambda(\epsilon,s),\dot{\Lambda}(\epsilon,s))}{\partial \dot{\Lambda}^i}\Big)\frac{d}{d\epsilon}\,
\Lambda^i(\epsilon,s)\,ds.
\label{firstvatiationtimearrivalfunctional}
\end{multline}
Note that in this formula $\epsilon$ has not been fixed to have the value $\epsilon=0$.
\section{\small{Second variation formula}}
\subsection{\small{The Chern connection on the pull-back bundle $\pi^*TM$}}
Our way to
introduce {\it Chern's connection} for Finsler spacetimes is as follows. The
 Cartan tensor components are defined by
\begin{align}
C_{ijk}:=\,\frac{1}{2} \frac{\partial g_{ij}}{\partial y^{k}}, \quad i,j,k=1,...,n,
\label{Cartantensor coefficients}
\end{align}
differently to the way it is introduced in \cite{BCS} by a factor $\frac{1}{F}$.
Therefore, because of the homogeneity of the tensor $g$, Euler's theorem implies
\begin{align}
C_{(x,y)}(y,\cdot,\cdot)=\,\,\frac{1}{2} y^k\,\frac{\partial g_{ij}}{\partial y^{k}}=0.
\label{homogeneityofcartan tensor}
\end{align}
The formal second kind Christoffel symbols
${\gamma}^{i}\,_{jk}(x,y)$ are defined by the expression
\begin{align*}
 {\gamma}^{i}\,_{jk}=\frac{1}{2}g^{is}(\frac{\partial g_{sj}}{\partial
x^{k}}-\frac{\partial g_{jk}}{\partial x^{s}}+\frac{\partial
g_{sk}}{\partial x^{j}}),\quad i,j,k,s=1,...,n.
\end{align*}
The {\it non-linear connection coefficients} are defined on ${ N}:=TM\setminus\{0\}$ to be
\begin{align*}
{N^{i}\,_{j}}={\gamma}^{i}\,_{jk}{y^{k}}-C^{i}\,_{jk}
{\gamma}^{k}\,_{rs}\,{y^{r}}{y^{s}},\quad
i,j,k,r,s=1,...,n,
\end{align*}
where $C^{i}\,_{jk}=g^{il}\,C_{ljk}$.

The connection coefficients define a connection on $N$ in the sense of Ehresmann, defining a decomposition
\begin{align}
TN=\mathcal{H}\oplus\mathcal{V},
\end{align}
where $\mathcal{V}=\,ker(d\pi)$.
An adapted frame to the above decomposition is determined by the smooth tangent basis for ${ T}_u { N},\, u\in { N}$:
\begin{align}
&&\left\{ \frac{{\delta}}{{\delta} x^{1}}|_u ,...,\frac{{\delta}}{{\delta} x^{n}} |_u, %%@
\frac{\partial}{\partial y^{1}} |_u,...,\frac{\partial}{\partial
y^{n}} |_u\right\},\,
 &&\frac{{\delta}}{{\delta} x^{j}}|_u =\frac{\partial}{\partial
x^{j}}|_u -N^{i}\,_{j}\frac{\partial}{\partial y^{i}}|_u ,\quad
i,j=1,...,n.
\label{adaptedframe}
\end{align}
Given a tangent vector $X\in\,T_x M$ and $u\in \pi^{-1}(x)$, there is a unique horizontal tangent
vector $h(X)\in \,T_u N$ with $d\pi(h(X))=X$ (horizontal lift of $X$).

The pull-back bundle $\pi^* {TM}$ is the maximal subsect of the cartesian product
${ N}\times{ TM}$ such that the diagram
\begin{align*}
\xymatrix{\pi^*{TM} \ar[d]_{\pi_1} \ar[r]^{\pi_2} &
{ TM} \ar[d]^{\pi_0}\\
{ N} \ar[r]^{\pi} & { M}}
\end{align*}
is commutative.
The projections on the first and second %%@
factors are
\begin{align*}
& \pi _1:\pi ^* { TM}\longrightarrow { N},\,
(u,\xi)\mapsto u,\quad \pi _2 :\pi ^* { TM}\longrightarrow { TM},\,
(u,\xi)\mapsto \xi.
\end{align*}
A {\it Chern type connection} is defined through the following
\begin{teorema}
Let $({ M},F)$ be a Finsler structure. The pull-back vector
bundle $\pi_1:{\pi}^{*}{TM}\rightarrow { N}$ admits a %%@
unique linear connection determined by the connection $1$-forms
 $\{{\omega}_j\,^i,\,\, i,j=1,...,n \} $ such that the following structure equations hold:
\begin{enumerate}
\item {\it Torsion free condition}
\begin{align}
d(dx^{i})-dx^{j}\wedge w_{j}\,^i=0,\quad i,j=1,...,n.
\label{torsionfree}
\end{align}
\item {\it Almost $g$-compatibility condition}
\begin{align}
dg_{ij}-g_{kj}w_{i}\,^k-g_{ik}w_{j}\,^k=2\,C_{ijk}{{\delta}y^{k}},\quad
i,j,k=1,...,n.
\label{almostgcompatibility}
\end{align}
\end{enumerate}
\label{chernconnection}
\end{teorema}
This result is proved along the same lines as in \cite{BCS} for the existence of
the Chern connection. This is because the hypothesis that one uses and the {\it Christoffel trick} is the same as in the positive case.
\begin{corolario}
Let $h(X)$ and $v(X)$ be the horizontal and vertical lifts of $X\in\,\Gamma TM$ to $TN$, and $\pi^*g$ the pull back-metric.
For the Chern connection the following properties hold:
\begin{enumerate}
\item  The almost $g$-compatibility metric condition is equivalent to
\begin{align}
{\nabla}_{v({{X}})}\pi^*g=2C({X},\cdot,\cdot),\quad
{\nabla}_{h({X})} \pi^*g=0,\quad {X}\in\,\Gamma TN.
\label{covariantalmostmetriccondition}
\end{align}

\item The torsion-free
condition of the
Chern connection is equivalent to the %%@
following:
\begin{enumerate}
\item Null vertical covariant derivative of sections of ${\pi}^*{TM}$:
\begin{align}
{\nabla}_{{v({X})}} {\pi}^* Y=0,
\label{covariantensorfreecondition1}
\end{align}
for any vertical component $v(X)$ of $X$.
\item  Let us consider $X,Y\in {TM}$ and their horizontal
lifts $h(X)$ and $h(Y)$. Then
\begin{align}
\nabla_{h(X)} {\pi}^* Y-{\nabla}_{h(Y)}{\pi}^*X-{\pi}^* ([X,Y])=0.
\label{covariantensorfreecondition2}
\end{align}
\end{enumerate}
\end{enumerate}
\end{corolario}
The proof of this corollary is similar to the positive case that one can be find in \cite{RF}.

 Given a Finsler structure, the curvature endomorphism of the Chern connection have two pieces different from
 zero \cite{BCS}. One of the pieces correspond to the $hh$-curvature and its components are given by
\begin{align}
& R^{i}\,_{jkl}  =\frac{\delta \Gamma^i\,_{jk}}{\delta x^l}-\frac{\delta \Gamma^i\,_{jl}}{\delta x^k}-\
 \Gamma^i\,_{hl}\Gamma^h\,_{jk}+\Gamma^i\,_{hk}\Gamma^h\,_{jl},\quad h,i,j,k,l=1,...,n.
\end{align}
\begin{lema}
Let $\lambda:[0,1]\longrightarrow M$ be a geodesic of $L$.
Then along  $\lambda:[0,1]\longrightarrow M$ the following relations are true for the Chern connection,
\begin{enumerate}
\item The connection preserves $g$ along $\lambda$ in the sense that
\begin{align}
\frac{d}{d r}\big(g_{\dot{\lambda}}(Y,W)\big)=\,g_{\dot{\lambda}}(\,\nabla_{h(X)}(\dot{\lambda})\,\pi^*Y,\pi^*W)+
\,\,g_{\dot{\lambda}}(\pi^*Y,
\,\nabla_{h(X)}(\dot{\lambda})\pi^*W),
\end{align}
holds for all $X,Y,W\in \Gamma TM.$
\item The connection $\nabla$ is torsion free: if $[X,Y]=0$, then
\begin{align}
\nabla_{h(X)}(\lambda,\dot{\lambda})\,\pi^*Y=\,\nabla_{h(Y)}(\lambda,\dot{\lambda})\,\pi^*X,\quad X\in \Gamma TM.
\end{align}
\end{enumerate}
\label{sobreprimeravaraicion}
\end{lema}
{\bf Proof.} The first property is the evaluation of the almost metric compatibility condition of the
Chern connection along $u(s)=(\lambda(s),\dot{\lambda})$, and with $\frac{d}{dr}$ the
derivative along the vector field $X\in \Gamma \pi^*TM$. When the condition (\ref{almostgcompatibility}) is evaluated at $u$, the left hand side is
\begin{align*}
(dg_{ij}-g_{kj}w_{i}\,^k-g_{ik}w_{j}\,^k)(X) & =dg_{ij}|_u(X)\,-g_{kj}|_u\,w_{i}\,^k(X)\,-g_{ik}|_u\,w_{j}\,^k(X)\\
& =\frac{d}{dr}(g_{ij}|_u\,(X))\,-g_{kj}|_u\,w_{i}\,^k(X)\,-g_{ik}|_u\,w_{j}\,^k(X).
\end{align*}
This expression is proportional to the Cartan tensor evaluated along $\lambda$ and contracted with $\dot{\lambda}$. along the causal geodesic $\lambda$ evaluated in the first entry at $\dot{\lambda}(s)$.
However, the Cartan tensor along a critical point of $L=-F^2$ is zero,
 \begin{align*}
 2\,y^k\,C_{ijk}(x,y)=0
\end{align*}
by homogeneity of degree zero of $g$ and by Euler's theorem. The second statement is a consequence of
the torsion free condition of the Chern connection.\hfill$\Box$
\begin{lema}
Let $\nabla$ be the Chern connection on $\pi^*TM\longrightarrow N$. The first variation of the functional $E$ is given by the expression
\begin{align}
 0 =\,g_{\dot{\Lambda}}(\dot{\Lambda}(1,\epsilon),{\gamma}'(\tau(\Lambda)))\,\frac{d}{d\epsilon}\,\big(\tau\Lambda(\epsilon,s)\big)\,
 -\int^1_0\,g_{\dot{\Lambda}}(\,\nabla_{h(\dot{\Lambda})}(\dot{\Lambda})\,\dot{\Lambda},\frac{d}{d\epsilon}\,\Lambda(\epsilon,s)\big)\,ds,
 \label{primeravariacionformula}
\end{align}
where $h(\dot{\Lambda})\in\,TN$ is the horizontal lift of $\dot{\Lambda}\in T_{\lambda}M$.
\label{primeravariacion}
\end{lema}
{\bf Proof}. If a
connection is torsion-free and preserves the metric along a geodesic (as it is the case
because of (\ref{covariantalmostmetriccondition}) and {\it lemma}  \ref{sobreprimeravaraicion}, one obtains the relation
\begin{align*}
\int^1_0\,ds\Big(\frac{\partial L(\Lambda(\epsilon,s),\dot{\Lambda}(\epsilon,s))}{\partial \Lambda^i}
\,&-\frac{d}{ds}\,\frac{\partial L(\Lambda(\epsilon,s),\dot{\Lambda}(\epsilon,s))}{\partial \dot{\Lambda}^i}\Big)\,\Lambda^i(s)=\,\int^1_0\,ds\frac{1}{2}
\,g_{\dot{\Lambda}}(\nabla_{h(\dot{\Lambda})}(\dot{\Lambda})\,\dot{\Lambda},\frac{d}{d\epsilon}\,\Lambda(\epsilon,s)),
\end{align*}
from which follows the result.\hfill $\Box$
\begin{definicion}
 Let $\lambda\in \,\mathcal{C}_{q,\gamma, c}$ and $\Lambda_1$, $\Lambda_2$ be
 two allowed variations of $\lambda$ such the corresponding vector fields along $\lambda$ are
  $A(s)=\,\frac{d\,\Lambda_1(s,\epsilon)}{d \epsilon}|_{\epsilon =0}$ and $B(s)=\,\frac{d\,\Lambda_2(s,\epsilon)}{d\,\epsilon}|_{\epsilon =0}$.
  We define the index form acting on $A(s)$ and $B(s)$ to be
\begin{align}
J_{\lambda}(A,B):=\,\int^1_0\,\Big(g_{\dot{\lambda}}\big(B(s),R_{\dot{\lambda}}(A(s),\dot{\lambda}(s))
(\dot{\lambda}(s))\big)-\,g_{\dot{\lambda}}(\nabla_{\dot{\lambda}}{A}(s),\nabla_{\dot{\lambda}}{B}(s))\Big)\,ds,
\label{formadelindice}
\end{align}
where $R$ is the $hh$-curvature of the covariant derivative $^+\nabla$ along $\lambda$  induced from the Chern connection.
\end{definicion}
Since the Cartan tensor along a geodesic is zero, the $hh$-curvature tensor
along $\lambda$ is equal to the Riemann curvature tensor along $\lambda$,
\begin{align}
R^{i}\,_{jkl}  =\frac{\partial \gamma^i\,_{jk}}{\partial x^l}-\frac{\partial \gamma^i\,_{jl}}{\partial x^k}-\,
 \gamma^i\,_{hl}\gamma^h\,_{jk}+\gamma^i\,_{hk}\gamma^h\,_{jl},\quad h,i,j,k,l=1,...,n.
\label{Riemanncurvaturetensoralonggeodesic}
\end{align}

One can relate the space
  $T_{\lambda}\mathcal{C}_{q,\gamma,c}$ with the space of tangent vector fields along the geodesic $\lambda$,
  \begin{align*}
\mathcal{T}_{\lambda}:=\Big\{\,& A:[0,1]\longrightarrow T M, \quad C^{\infty} \,\textrm{functions}\, s.t.\\
& 1. \pi(A)=\lambda,\\
& 2. \,A(0)=0,\\
& 3. \,A(1)=0,\\
& 4. \,g_{\dot{\lambda}}\,(\dot{\lambda}(s),A(s))=0\Big\}.
\end{align*}
 \begin{lema}
 Given a causal geodesic $\lambda:[0,1]\longrightarrow \,M$,
 the vector spaces  $\mathcal{T}_{\lambda}$ and $T_{\lambda}\mathcal{C}_{q,\gamma,c}$ are isomorphic.
 \label{lemaquerelacionalosindices}
 \end{lema}
{\bf Proof}. We prove that the vector spaces coincide. For both spaces $A(0)=0$.
 From property $2$ in the definition of $\mathcal{C}_{q,\gamma,c}$ one has that for
any $A\in\, T_{\lambda}\mathcal{C}_{q,\gamma,c}$ and homogeneity of $L$
\begin{align*}
\frac{d}{ds}L\Big(\Lambda(s,\epsilon),\dot{\Lambda}(s,\epsilon)\Big)\Big|_{\epsilon=0}  =\,&\frac{\partial \,L}{\partial \lambda^i(s)}\,A^i(s)\,
+\frac{\partial \,L}{\partial \dot{\lambda}^i(s)}\,\dot{A}^i(s)\,+g_{ij}(\lambda,\dot{\lambda}(s))\dot{\lambda}^j(s)\,\dot{A}^i(s)=0.
\end{align*}
Using the condition that $\lambda$ is a geodesic this reduces to
\begin{align*}
\frac{d}{ds}L\Big(\Lambda(s,\epsilon),\dot{\Lambda}(s,\epsilon)\Big)\Big|_{\epsilon=0}
=\,&\frac{d}{ds}\Big(\frac{\partial \,L}{\partial \dot{\lambda}^i(s)}\Big)\dot{A}^i
+\,g_{ij}(\lambda,\dot{\lambda}(s))\dot{\lambda}^j(s)\,\dot{A}^i(s)\\
& =\, \frac{d}{ds}\Big(g_{ij}(\lambda,\dot{\lambda}(s))\dot{\lambda}^j(s)\,{A}^i(s)\Big)=0.
\end{align*}
This {\it conservation law} proves that $A(s)$ satisfies point $3$. To check point $3$, let us note that $A(1)= a\gamma'(\tau(\lambda))$.
By point $3$,
\begin{align*}
0=\,g_{\dot{\lambda}}(A(1),\dot{\lambda}(1))=\,g_{\dot{\lambda}}(a\frac{d}{dt}\gamma(\tau(\lambda)), \dot{\lambda}(1))
=\,a \,g_{\dot{\lambda}}(\gamma'(\tau(\lambda)), \dot{\lambda}(1)).
\end{align*}
 In the last expression the
factor $\,g_{\dot{\lambda}}(\gamma'(\tau(\lambda)), \dot{\lambda}(1))$ cannot be zero,
since $T$ and $\dot{\gamma}$ are both not mutually non-orthogonal timelike vectors and the
hypothesis that 
\begin{align*}
g_{\dot{\lambda}}(\dot{\lambda},T(\lambda(1)))<0.
\end{align*}
 Therefore, $a=0$ and the result is proved.\hfill$\Box$

{\it Lemma} \ref{lemaquerelacionalosindices} is used in the proof of the second variation formula and in the index formula in the next {\it section}.
\subsection{\small{Second variation of the time arrival functional}}
Note that for each allowed variation $\Lambda(\epsilon,s)$ the time arrival functional is
a smooth function on the variable $\epsilon$. Therefore, the hessian is defined along a critical point.
\begin{proposicion}
Let $\lambda:[0,1]\longrightarrow M$ be in $\mathcal{C}_{q,\gamma, c}$ a causal geodesic.
Let $\Lambda(\epsilon,s)$ an allowed variation with $A=\frac{d}{d\epsilon}\big|_{\epsilon=0}\Lambda(\epsilon,s)$. Then
\begin{align}
\frac{d^2}{d\epsilon^2}\,\tau(\Lambda(\epsilon,\cdot))\Big|_{\epsilon=0}\,=\,\frac{J_{\lambda}(A,A)}
{g_{\dot{\lambda}}({\gamma}'(\tau(\lambda)),\dot{\lambda}(1))}.
\label{segundavariacionformula}
\end{align}
\label{segundavariacion}
\end{proposicion}
{\bf Proof.} From the first variation formula (\ref{primeravariacionformula}) we take the second derivative respect to $\epsilon$,
\begin{multline*}
0  =\,\frac{d}{d\,\epsilon}\,\Big(g_{\dot{\Lambda}}(\dot{\Lambda}(\epsilon,1),{\gamma}'(\tau(\Lambda)))\,
\frac{d}{d\epsilon}\,\big(\tau(\Lambda(\epsilon,s))\big) -\int^1_0\,2\,\,g_{\dot{\Lambda}}
(\,\nabla_{h(\dot{\Lambda})}(\dot{\Lambda})\,\Lambda,\frac{\partial}{\partial \epsilon}\,
\dot{\Lambda}(s,\epsilon)\big)\,ds\Big)\Big|_{\epsilon =0}\\
=\,\frac{d}{d\,\epsilon}\,\Big(g_{\dot{\Lambda}}(\dot{\Lambda}(\epsilon,1),{\gamma}'(\tau(\Lambda)))\Big)\Big|_{\epsilon =0}
\,\frac{d}{d\epsilon}\Big|_{\epsilon =0}\,\big(\tau\Lambda(\epsilon,s)\big)
 +\,g_{\dot{\Lambda}}(\dot{\Lambda}(\epsilon,1),{\gamma}'(\tau(\Lambda)))\frac{d}{d\,\epsilon}\Big|_{\epsilon =0}\,
\Big(\frac{d}{d\epsilon}\,\big(\tau\Lambda(\epsilon,s)\big)\Big)
\\- \,\frac{d}{d\,\epsilon}\,\Big(\int^1_0\,g_{\dot{\Lambda}}(\,\nabla_{h(\dot{\Lambda})}(\dot{\Lambda})\,\dot{\Lambda},
\frac{\partial}{\partial \epsilon}\,\Lambda(\epsilon,s)\big)\,ds\Big)\Big|_{\epsilon =0}.
\end{multline*}
Evaluated on a geodesic, the first term is zero, since by Proposition~\ref{fermatprinciplefortimelikegeodesics}
\[\frac{d}{d\epsilon}\big|_{\epsilon=0}\,\big(\tau(\Lambda(\epsilon,s))\big)=0.\]
Using the metric compatibility and the torsion free
conditions along $\lambda$, the condition $[\dot{\Lambda},\frac{\partial}{\partial \epsilon }{\Lambda}]=0$
and the geodesic equation in terms of the connection
$\nabla_{\dot{\lambda}}(\dot{\lambda})\,\dot{\lambda}=0$, one obtains the following expression
\begin{align*}
\frac{d}{d\,\epsilon}\,\Big(\,&\int^1_0\,g_{\dot{\Lambda}}(\,\nabla_{h(\dot{\Lambda})}(\dot{\Lambda})\,\dot{\Lambda},
\frac{d}{d\epsilon}\,
\Lambda(\epsilon,s)\big)ds\Big)\Big|_{\epsilon =0}\\&=\,\int^1_0\,\frac{d}{d\,\epsilon}\,\Big(\,g_{\dot{\Lambda}}
(\,\nabla_{h(\dot{\Lambda})}(\dot{\Lambda})\,\dot{\Lambda},
\frac{d}{d\epsilon}\,\Lambda(s)\big)ds\Big)\Big|_{\epsilon =0}\\
& =\int^1_0\,\,g_{\dot{\lambda}}\big(\,\nabla_{h({A})}(\dot{\lambda})\,\nabla_{h(\dot{\lambda})}(\dot{\lambda})\,\dot{\lambda},A\big)ds
 +\,\int^1_0\,\,g_{\dot{\lambda}}\big(\,\nabla_{h(\dot{\lambda})}(\dot{\lambda})\,\dot{\lambda},\,\nabla_{h({A})}(\dot{\lambda})\,A\big)ds\\
& =\,\int^1_0\,\,g_{\dot{\lambda}}\big(\,\nabla_{(A)}(\dot{\lambda})\,\nabla_{h(\dot{\lambda})}(\dot{\lambda})\,\dot{\lambda},A\big)ds\\
& =\,\int^1_0\,\,g_{\dot{\lambda}}\big(\,\nabla_{h(\dot{\lambda})}(\dot{\lambda})\,\nabla_{h(\dot{A})}(\dot{\lambda})\,\dot{\lambda}+
\,R(A,\dot{\lambda})\dot{\lambda},A\big).
\end{align*}
In the last line, the first term can be computed more explicitly:
\begin{align*}
\int^1_0\,\,g_{\dot{\lambda}}\big(\,\nabla_{h(\dot{\lambda})}(\dot{\lambda})\,\nabla_{h({A})}(\dot{\lambda})\,\dot{\lambda},A\big)ds
&\,=\,-\int^1_0\,\,g_{\dot{\lambda}}\big(\,\nabla_{h({A})}(\dot{\lambda})\,\dot{\lambda},\,\nabla_{h(\dot{\lambda})}(\dot{\lambda})\,A\big)ds\\
& +\,\int^1_0\,\frac{d}{ds}\Big(g_{\dot{\lambda}}\big(\,\nabla_{h(\dot{\lambda})}(\dot{\lambda})\,A,A\big)\Big)ds\\
& \,=\,-\int^1_0\,\,g_{\dot{\lambda}}\big(\,\nabla_{h(\dot{\lambda})}(\dot{\lambda})\,A,\,\nabla_{h(\dot{\lambda})}(\dot{\lambda})\,A\big)ds,
\end{align*}
where in the last equality we have used Lemma~\ref{lemaquerelacionalosindices}.
Combining this last relation with the definition of the index $J_{\lambda}$ one obtains the result.\hfill $\Box$

The proof of {\it proposition} \ref{segundavariacion} suggests to compile the following properties in a {\it lemma},
\begin{lema}
Let $\lambda\in\mathcal{C}_{q,\gamma, c}$  be a causal geodesic. Then along any geodesic $\lambda$ the following properties hold:
\begin{enumerate}
\item The connection is {\it torsion-free} along $\lambda$.
\item The Cartan tensor is zero along a geodesic. Therefore, it is metric compatible along $\lambda$.
\item The $hh$-curvature reduces to the Riemann curvature (\ref{Riemanncurvaturetensoralonggeodesic}) along a geodesic.
\end{enumerate}
\label{Lorentzianproperties}
\end{lema}
{\bf Proof}. For Finslerian quantities, the base point vector is fixed along a geodesic.
Then
\begin{enumerate}
\item  For the Chern's connection, the connection is torsion-free on whole $N$.
\item  The Cartan tensor along a geodesic is zero \cite{BCS}.
\end{enumerate}
 As consequence that the Cartan tensor is zero,
the $hh$-curvature reduces to a formal Riemann curvature \ref{Riemanncurvaturetensoralonggeodesic}.\hfill$\Box$

\section{\small{Applications}}
 As a consequence of {\it lemma}
\ref{Lorentzianproperties}, most of the proofs from \cite[Chapter 10]{BEE}, \cite{Perlick} can be adapted to Finsler spacetimes with minimal changes.
Using this fact, we provide two applications. First, sice the relation between the second variation
of time functional and the index form along a causal geodesic is the same
than in \cite{Perlick}, one can study the character of the critical points of the time
arrival functional as in the Lorentzian case. These results serve to illustrate that indeed one can {\it transplant} the
methods of the Lorentzian case to the Finsler spacetime category.
Also note that the result is valid for causal geodesics, not only timelike geodesics.
Second, one can translate the techniques from \cite{BEE} to obtain a Morse
index theorem for timelike geodesics of time orientable Finsler spacetimes. This is related with the
index of the Hessian of the time arrival functional.
\subsection{\small{The character of the critical points of the time arrival functional}}
Let us consider the vector spaces
\begin{align*}
& V^{\bot}(\lambda):=\,\{\textrm{
all piecewise smooth vector fields $A$ along $\lambda$
\,s.t.\, $g_{\dot{\lambda}}(\dot{\lambda},A(s))=0$}\}\\
& V^{\bot}_0(\lambda):=\{A\in\,V^{\bot},\,s.t. \, A(0)=A(1)=0\}.
\end{align*}
 A direct application
of {\it lemma} \ref{Lorentzianproperties} is that for arbitrary smooth vector fields along $\lambda$, the index form (\ref{formadelindice}) is
\begin{multline}
J_{\lambda}(A,B)=\\-g_{\dot{\lambda}}(\nabla_{h(\dot{\lambda})}{A},B)|^1_0\,+,\int^1_0\,\Big(g_{\dot{\lambda}}\big(B(s),
\nabla_{h(\dot{\lambda})}\nabla_{h(\dot{\lambda})}A(s)+R_{\dot{\lambda}}(A(s),\dot{\lambda}(s))
(\dot{\lambda}(s))\big)\Big)\,ds.
\label{formadelindice2}
\end{multline}
If $B\in \, V^{\bot}_0(\lambda)$ and $A$ are smooth vector fields along $\lambda$, the index form is
\begin{align}
J_{\lambda}(A,B)=\int^1_0\,\Big(g_{h(\dot{\lambda})}\big(B(s),
\nabla_{h(\dot{\lambda})}\nabla_{h(\dot{\lambda})}A(s)+R_{\dot{\lambda}}(A(s),\dot{\lambda}(s))
(\dot{\lambda}(s))\big)\Big)\,ds,
\label{formadelindice3}
\end{align}
\begin{definicion}
Let $(M,L)$ be a Finsler spacetime and $\, \mathcal{C}_{q,\gamma,c}\ni\,\lambda:[a,b]\longrightarrow M$ a causal geodesic.
\begin{enumerate}
\item
A Jacobi field is a vector field $Y$ along $\lambda$ such that is a solution of the Jacobi equation
\begin{align}
\nabla_{h(\dot{\lambda})}\,\nabla_{h(\dot{\lambda})}\, Y+\,R_{\dot{\lambda}}(Y,\dot{\lambda})\dot{\lambda}=\,0,
\label{jabociequation}
\end{align}

\item Let $ \mathcal{C}_{q,\gamma,c}\ni\,\lambda:[a,b]\longrightarrow M$ be a causal geodesic. Then $\lambda(t_1)$ and $\lambda(t_2)$
are conjugate points along $\lambda$ iff there is a non-zero Jacobi field
such that $Y(t_1)=\,Y(t_2)=0$.
\label{jacobifieldconjugatepoints}
\end{enumerate}
\end{definicion}
\begin{lema}
Let $\lambda\in\mathcal{C}_{q,\gamma, c}$  be a causal geodesic without conjugate points. Then $J_{\lambda}(A,A)<0$ for any  $A\in\,V^{\bot}_0(\lambda)$.
\label{lemasobreelindice}
\end{lema}
{\bf Proof}. We learnt from the proof of {\it proposition} \ref{segundavariacion} and {\it lemma} \ref{Lorentzianproperties}
that Lorentz\-ian properties along geodesics carry over from the Lorentz\-ian case to Finsler spacetime case.
This makes the proof of the {\it lemma} completely analogous to the Lorentzian case (see \cite{BEE}
for the timelike case and \cite{BEE, Perlick} for the light-like case). \hfill$\Box$

The following proposition is a restatement of known results in Lorentzian geometry \cite{BEE},
\begin{proposicion}
Let $(M,L)$ be a Finsler spacetime and $\lambda:[0,1]\longrightarrow M\in\, \mathcal{C}_{q,\gamma,c}$ a causal geodesic and
$Y$ a Jacobi field along $\lambda$. Then

\begin{enumerate}
\item  The function along $\lambda$ given by $g_{\dot{\lambda}}(Y,\dot{\lambda})$ is an affine function.

\item Let $Y$ such that $Y(t_0)=Y(t_1)=0$ for different $t_1,t_2\in\,[0,1]$. Then
\begin{enumerate}
\item If $\lambda$ is a timelike geodesic, then $Y\in\,V^{\bot}_0(\lambda)$.

\item If $\lambda$ is a lightlike geodesic, then $Y$ is either orthogonal or parallel
to $\dot{\lambda}$ (therefore, it is lightlike).
\end{enumerate}
\item Let $Y$ such that $Y(0)=Y(1)=0$ for different $t_1,t_2\in\,[0,1]$. Then
\begin{enumerate}
\item If $\lambda$ is a timelike geodesic, then $\nabla_{h(\dot{\lambda})}\,Y\in\,V^{\bot}_0(\lambda)$.

\item If $\lambda$ is a lightlike geodesic, then $\nabla_{h(\dot{\lambda})}\,Y$ is orthogonal
to $\dot{\lambda}$ (therefore, it is lightlike).
\end{enumerate}
\end{enumerate}
\label{morelorentzianproperties}
\end{proposicion}
{\bf Proof}. The proof follows closely the proof for {\it lemma 10.9}, {\it corollary 10.10}
and {\it 10.11} in \cite{BEE}, through the use
of {\it lemma} \ref{Lorentzianproperties}.\hfill$\Box$

The following theorem is proved in a similar way as in \cite{Perlick},
\begin{teorema}
Let $\mathcal{C}_{q,\gamma,c}\ni\,\lambda:[a,b]\longrightarrow M$ be a causal geodesic and $\tau$ the time arrival functional. Then
\begin{enumerate}
\item If $\lambda$ does not have conjugate points, then it is a local minimum of $\tau$.

\item If $\lambda$ has intermediate conjugate points, then it is a local saddle point of $\tau$.
\end{enumerate}
\label{characteroftheconjugatepoint}
\end{teorema}
{\bf Proof}. The proof of the first statement follows the same steps than the Lorentzian case (see \cite{Perlick}): from the formula for the second
variation of the arrival time functional (\ref{segundavariacionformula}),
 it follows that if $J_{\lambda}<0(A;A)$ (by {\it lemma} \ref{lemasobreelindice}), then the time arrival is a local minimal .

 The proof of the second statement is identical to the proof in \cite{BEE} for the case of
 timelike geodesics and to \cite{Perlick} for lightlike geodesics and
 will not be rewrite here. Note that it is essential in the proof both {\it lemma}
 \ref{Lorentzianproperties} and {\it proposition} \ref{morelorentzianproperties},
 which provide exactly the same tools as in the Lorentzian case.\hfill$\Box$

 \subsection{\small{Morse index theorem for the arrival time functional in Finsler spacetime spacetimes}}
 The common feature of the results and techniques in the Finsler and Lorentzian case suggests
  that there is also a Finsler spacetime version of the Morse index theorem. Indeed,
  there is such result as we indicate below.
The index of $\lambda$ is equal to the number on conjugate points along $c$ counted
with multiplicity, that is, counting the dimension of the vector space of
 Jacobi fields $J_i\in\,V^{\bot}_0(\lambda)$ vanishing at each conjugate point $\lambda(s_i)$,
\begin{align}
I(\lambda):=\sum_{s\in\,(0,1)} dim (J_i).
\end{align}
The index of the bilinear form $\tau$ along $\lambda$ denoted by $I(\tau,\lambda)$ is the
supreme of the dimensions of all subspaces of $V^{\bot}_0(\lambda)$ on which the Hessian $Hess(\tau)$ is negative. Then one has the following result,
\begin{teorema}
Given a timelike geodesic $\lambda:[0,1]\longrightarrow M$, the number of conjugate points $I(\lambda)$ is given by
\begin{align}
I(\lambda)=\,I(\tau,\lambda)=\sum_{s\in\,(0,1)} dim (J_i).
\label{index}
\end{align}
\label{teoremadelindice}
\end{teorema}
{\bf Proof}. The methods and results described before suggest the existence of a version of the Morse index theorem for the functional energy
for general Finsler spacetime . That this is the case will be seen in \cite{GallegoTorromePiccioneVictorio}.
By {\it lemma} \ref{lemaquerelacionalosindices}, one relates such index with $I(\tau,\lambda)$ by formula (\ref{segundavariacionformula}).\hfill$\Box$

\begin{comentario}
 It is worth to mention that one can establish {\it theorem } \ref{characteroftheconjugatepoint}
 as a direct consequence of the index theorem \ref{teoremadelindice}.
\end{comentario}

\section{Discussion}

In this paper we have investigated  a generalization of Fermat's principle for causal curves in time oriented
Finsler spacetimes and some related results. The second
variation formula and Morse index theorem associated with the time arrival functional have also been investigated. One notes easily that there is not formal difference with the corresponding Lorentzian results \cite{Perlick}. However, this formal analogy is only apparent, since the geometric objects in Finsler spacetime lives on $N$ and not on $M$ directly. Still, 
we have shown that Beem's theory of Finsler spacetimes offers an adequate framework to generalize the standard theory of Lorentzian geometry to Finsler spacetimes.

There are several questions considered in this paper that deserve further discussion. For instance, the framework discussed here is 
 a restriction of a more general formulation of Finsler spacetimes. In particular,
 the question of the choice of the base point where the  geometric quantities are evaluated
 is non-trivial (for a discussion of this issue see \cite{Rund} and \cite{Ishikawa}). When we define
the notion of causal curves $\lambda:I\to M$, the base-point $(x,y)\in\,TM\setminus\{0\}$ where the function $L(x,y)$ is evaluated is
\begin{align*}
(\lambda(s),y(s))=(\lambda(s),\dot{\lambda}(s)),\,\dot{\lambda}\in\,T_{\lambda}M.
\end{align*}
This choice is not necessary, as it was pointed out by Ishikawa
  \cite{Ishikawa}. On the other hand, it is a natural choice, since the expressions appearing in the formulation
 of Fermat's principle, the second variation formula and the Morse
  index theorem must be evaluated along the curves $\lambda:I\to M$ at $(\lambda(s),\dot{\lambda}(s))$.

It is also convenient to use weaker regularity conditions for the Lagrangian $L$ than the required in Beem's definition. This is for instance the situation in
the examples 2.10 to 2.12, where although motivated from physical models, they are not regular on all the slit tangent bundle $N=TM\setminus\{0\}$.
Example 2.12 is indeed quite pathological,
since in this case the Lagrangian $L$ is singular in the full null cone. A convenient theory to deal with some of such singular examples is the theory of conic metrics \cite{JavaloyesSanchez}. However, even such general framework is not enough to deal with all the interesting problems appearing from physical applications.

Of particular importance has been the fact that the Finsler spacetime is time orientable. At present, this hypothesis is necessary in the proof of Fermat's principle \ref{fermatprinciplefortimelikegeodesics}. One wonders if there is a proof of Fermat's principle without the requirement of time orientability. Indeed, such proof exists for lightlike Finsler geodesics \cite{Perlick06}. This question for the timelike case remains open. 

Similarly, one can investigate the Morse index theorem for timelike geodesics of a Finsler spacetime in complete analogy to the Lorentzian case. The extension to lightlike geodesics must be done more carefully, because in the case of Finsler spacetimes, additional singularities could appear in the light cones.

 As we said before, there is a hierarchy of {\it reversible conditions}. Our definition of reversible Lagrangian \ref{definitionreversivilityofL} is stronger than the reversibility condition in Beem's theory, which is still stronger than the definition of reversibility in Pfeifer-Wohlfart's theory \cite{PfeiferWohlfarth}. The merit of our definition is based on the Lagrangian $L$, that is the function that appears in the formulation of Fermat's principle and that we emphasize as the fundamental object (together with $M$) in the definition of a Finsler spacetime and that is applicable to any causal curve or vector (in contrast with Beem's definition, that relies on the Finsler function $F(x,y)$, which is non-singular on the lightcone). The hierarchy described before also justifies why one does need to restrict to reversible spacetimes, since it could be a strong assumption, eliminating interesting spaces by being non-reversible.
 
 We did not restrict our considerations to reversible Finsler spacetimes in the sense of definition \ref{definitionreversivilityofL}. However, we have in mind that non-reversibility of the Lagrangian could produce additional difficulties for a consistent causal structure of $L$. In particular, one should impose conditions preventing possible almost-everywhere smooth, continuous causal loops.

 Let us mention that the Fermat's principle for Finsler spacetimes, the second variation formula and the Morse index theorem,
 apply not only in Finsler spacetimes describing gravity, but also to geometric models that
 describe the motion of particles under the action of gravity in combination with other interactions. This type of models include the motion in locally anisotropic media. 

 \subsection*{Acknowledgements} We would like to acknowledge to M. A. Javaloyes and V. Perlick for their useful comments and interest on this work.

\small{
}

\begin{thebibliography}{22}
\bibitem{Asa} G. S. Asanov, {\it Finsler
Geometry, Relativity and Gauge Theories}, D. Reidel, Dordrecht.
(1985).
\bibitem{BCS} D. Bao, S. S. Chern and Z. Shen, {\it An Introduction to
Riemann-Finsler Geometry}, Graduate Texts in Mathematics 200,
Springer-Verlag (2000).
\bibitem{Beem1} J. K. Beem, {\it Indefinite Finsler Spaces and Timelike Spaces},
Canad. J. Math. {\bf 22}, 1035 (1970).
\bibitem{BEE} John K. Beem, P. E. Ehrlich, K. L. Easly,
{\it Global Lorentzian Geometry}, Second Edition, CRC Press (1996).
\bibitem{Bogoslosvky} G. Bogoslosvky, {\it Some physical displays of
the space anisotropy relevant to the feasibility of its being detected at laboratory}, arXiv:0706.2621.
\bibitem{Chicone} C. Chicone, {\it Ordinary
Differential Equations with Applications}, 2nd. Edition, Springer-Verlag (2006).
\bibitem{CohenGlashow} A. G. Cohen, S. L. Glashow, {\it Very special relativity}, Phys. Rev. Lett. {\bf 97}, 021601 (2006).
\bibitem{GallegoTorrome} R. Gallego Torrom\'e, {\it
    Averaged Dynamics Associated with the Lorentz Force Equation}, arXiv:0905.2060;
{\it Averaged dynamics of ultra-relativistic charged particles beams}, PhD. Thesis, Lancaster University (2010).
\bibitem{RF} R. Gallego Torrom\'e, F. Etayo, {\it On a rigidity condition for Berwald spaces}, RACSAM 104 (1), 69-80 (2010).
\bibitem{GallegoTorromePiccioneVictorio} R. Gallego Torrom\'e, P. Piccione, H. Vit\'orio, work in progress.
\bibitem{GiannoniMasielloPiccione} F. Giannoni, A. Masiello, P. Piccione, {\it A timelike extension of Fermat's
 principle in general relativity and applications}, Cal. Var. 6, 263-283 (1998).
 \bibitem{GibbonsGomisPope} G. W. Gibbons, J. Gomis, C. N. Pope, {\it General very special relativity is Finsler geometry}, Phys.Rev.D76:081701, (2007).
\bibitem{GirelliLiberatiSindoni} F. Girelli, S. Liberati and L. Sindoni, {\it Phys. Rev. D} {\bf 75} 0604015 (2006).
\bibitem{HawEllis} S. W. Hawking and G. F. R. Ellis,
{\it The Large Scale Structure of the Space-Time},
Cambridge Monographs on Mathematical Physics (1973).
\bibitem{Ishikawa} H. Ishikawa, {\it Note on Finslerian relativity}, J. Math. Phys. {\bf 22}, 995 (1981).
\bibitem{JavaloyesSanchez} M. A. Javaloyes, M. Sanchez, {\it On the definition and examples of Finsler metrics}, arXiv: 1111.5066.
\bibitem{Kovner} I. Kovner, {\it Fermat principle in arbitrary gravitational fields}, Astrophysical Journal, vol. 351, March 1, p. 114-120 (1990).
\bibitem{Kosteleky} A. Kosteleky, {\it Riemann-Finsler geometry and Lorentz-violating kinematics}, 	 Phys.Lett.B701:137-143 (2011).
\bibitem{KouretsisStathakopoulosStraviros} A. P. Kouretsis, M. Stathakopoulos and P. C. Straviros, {\it Imperfect fluids, Lorentz violations and Finsler cosmology}, Phys. Rev. D {\bf 82} (2010) 064035.
\bibitem{LämmerzahlLorekDittus} C. L$\ddot{a}$mmerzahl, D. Lorek, H. Dittus,
{\it Confronting Finsler spacetime with experiment}, Gen. Rel. Grav. {\bf 41}, 1345-1353 (2009).
\bibitem{MironAnastasiei} R. Miron and M. Anastasiei, The Geometry of Lagrange Spaces: Theory and Applications,
Vol. 59, Fundamental Theories of Physics, Kluwer Academic Publishers (1994).
\bibitem{Perlick} V. Perlick, {\it On Fermat's principle in general relativity: I. The general case}, Class. Quantum Grav. 7, 1319-1331 (1990).
\bibitem{Perlick2000} V. Perlick, {\it Ray Optics, Fermat's principle and applications to general relativity}, Springer Heidelberg (2000).
 \bibitem{Perlickliving} V. Perlick, {\it Gravitational lensing from a spacetime perspective},
 Living Rev.\ Relativity 7 (2004), 9. URL :http://www.livingreviews.org/lrr-2004-9/.
\bibitem{Perlick06} V. Perlick, {\it Fermat Principle in Finsler Spacetimes},
 Gen. Rel. Grav. {\bf 38} 365-380  (2006).
\bibitem{PfeiferWohlfarth} C. Pfeifer and M. N. R. Wohlfarth, {\it Causal structure and electrodynamics on Finsler spacetimes},
Phys. Rev. D {\bf 84}:044039 (2011).
\bibitem{Randers} G. Randers, {\it On an Asymmetrical Metric in the
Four-Space of General Relativity}, Phys. Rev. {\bf 59},
195-199 (1941).
\bibitem{Rund} H. Rund, {\it The geometry of Finsler spaces}, Die Grundlehren der Mathematischen Wissenchaften, Band 101, Springer Verlag (1959).
\bibitem{Rutz} S. F. Rutz, {\it A Finsler Generalization of Einstein's Vacuum Field Equation}, Gen. Rel. Grav., Vol. 25, nº 11, 1139-1158 (1993).
\bibitem{SkakalaVisser} J. Skakala, M. Visser,{\it Pseudo-Finslerian spacetimes and multi-refrigence},
International Journal of Modern Physics D{\bf 19}, 1119-1146 (2010).
\bibitem{SkakalaVisser2} J. Skakala, M. Visser, {\it Bi-metric pseudo-Finslerian spacetimes}, Journ. Geom. Phys. {\bf 61}, 1386-1400 (2011).
 \bibitem{Warner} F. Warner, {\it Foundations of Differentiable Manifolds and Lie Groups}, Scott, Foresman and Company (1971).
 \bibitem{Weyl} H. Weyl, Ann. Phys. Lpz. {\bf 54} 117 (1917).
\end{thebibliography}
\end{document}